\documentclass[a4paper]{amsart}  

\usepackage{amssymb} \usepackage{amscd} \usepackage{times}

\def\zbb{\mathbb{Z}}  
  
  \def\phi{\varphi}
 \def\p1{{\mathbb{P}^1_\zbb}}

\newtheorem{Theorem}{\quad Theorem}[section]

\newtheorem{Lemma}[Theorem]{\quad Lemma}

\begin{document}

\title{ A uniform estimate for an equation with Holderian condition and boundary singularity.}

\author{Samy Skander Bahoura}

\address{Departement de Mathematiques, Universite Pierre et Marie Curie, 2 place Jussieu, 75005, Paris, France.}
              
\email{samybahoura@yahoo.fr, samybahoura@gmail.com} 

\date{}

\maketitle

\begin{abstract}

We consider the following problem on open set  $ \Omega $ of $ {\mathbb R}^2 $:

\begin{displaymath}  \left \{ \begin {split} 
      -\Delta u_i & =  |x-x_0|^{-2\alpha}V_i e^{u_i} \,\, &&\text{in} \!\!&&\Omega \subset {\mathbb R}^2, \\
                  u_i  & = 0  \,\,             && \text{in} \!\!&&\partial \Omega.               
\end {split}\right.
\end{displaymath}

Here, $ x_0 \in \partial \Omega $ and, $ \alpha  \in (0, 1/2) $. 

\bigskip

We assume, for example that:

$$ \int_{\Omega}|x-x_0|^{-2\alpha}V_i e^{u_i} dy \leq 16\pi-\epsilon , \,\, \epsilon >0 $$

1) We give, a quantization analysis of the previous problem under the conditions:

$$ \int_{\Omega}  |x-x_0|^{-2\alpha}e^{u_i} dy  \leq C, $$

and,

$$ 0 \leq V_i \leq b  < + \infty $$

2) In addition to the previous hypothesis we assume that $ V_i  $ $ s- $ holderian with $ 1/2 < s \leq 1$, then we have a compactness result, namely:

$$ \sup_{\Omega} u_i \leq c=c(b, C, A, s, \epsilon, \alpha,x_0, \Omega). $$

where $ A $ is the holderian constant of $ V_i $.

\end{abstract}

\section{Introduction and Main Results} 

We set $ \Delta = \partial_{11} + \partial_{22} $  on open set $ \Omega $ of $ {\mathbb R}^2 $ with a smooth boundary.

\bigskip

We consider the following problem on $ \Omega \subset {\mathbb R}^2 $:

\begin{displaymath} (P) \left \{ \begin {split} 
      -\Delta u_i & =   |x-x_0|^{-2\alpha}V_i e^{u_i} \,\, &&\text{in} \!\!&&\Omega \subset {\mathbb R}^2, \\
                  u_i  & = 0  \,\,             && \text{in} \!\!&&\partial \Omega.               
\end {split}\right.
\end{displaymath}

Here, $ x_0 \in \partial \Omega $ and, $ \alpha  \in (0, 1/2) $.  

\bigskip

We assume that,

$$ 0 \leq V_i \leq b  < + \infty, \,\, \int_{\Omega}|x-x_0|^{-2\alpha}e^{u_i} dy  \leq C, \,\, u_i \in W^{1,1}_0(\Omega) $$

The above equation is called, the Prescribed Scalar Curvature equation in
relation with conformal change of metrics. The function $ V_i $ is the
prescribed curvature.

\bigskip

Here, we try to find some a priori estimates for sequences of the
previous problem.

\bigskip

Equations of this type (in dimension 2 and higher dimensions) were studied by many authors, see [1-24]. We can
see in [8], different results for the solutions of those type of
equations with or without boundaries conditions and, with minimal
conditions on $ V $, for example we suppose $ V_i \geq 0 $ and  $ V_i
\in L^p(\Omega) $ or $ V_ie^{u_i} \in L^p(\Omega) $ with $ p \in [1,
+\infty] $. 

Among other results, we  can see in [8], the following important Theorem,

\bigskip

{\bf Theorem A}{\it (Brezis-Merle [8])}.{\it If $ (u_i)_i $ and $ (V_i)_i $ are two sequences of functions relatively to the previous problem $ (P) $ with, $ 0 < a \leq V_i \leq b < + \infty $, then, for all compact set $ K $ of $ \Omega $,

$$ \sup_K u_i \leq c = c(a, b, K, \Omega). $$}

A simple consequence of this theorem is that, if we assume $ u_i = 0 $ on $ \partial \Omega $ then, the sequence $ (u_i)_i $ is locally uniformly bounded. We can find in [8] an interior estimate if we assume $ a=0 $, but we need an assumption on the integral of $ e^{u_i} $. We have in [8]:

\smallskip

{\bf Theorem B} {\it (Brezis-Merle [8])}.{\it If $ (u_i)_i $ and $ (V_i)_i $ are two sequences of functions relatively to the previous problem $ (P) $ with, $ 0 \leq V_i \leq b < + \infty $, and,

$$ \int_{\Omega} e^{u_i} dy  \leq C, $$

then, for all compact set $ K $ of $ \Omega $,

$$ \sup_K u_i \leq c = c(b, C, K, \Omega). $$}

If, we assume $ V $ with more regularity, we can have another type of estimates, $ \sup + \inf $. It was proved, by Shafrir, see [23], that, if $ (u_i)_i, (V_i)_i $ are two sequences of functions solutions of the previous equation without assumption on the boundary and, $ 0 < a \leq V_i \leq b < + \infty $, then we have the following interior estimate:

$$ C\left (\dfrac{a}{b} \right ) \sup_K u_i + \inf_{\Omega} u_i \leq c=c(a, b, K, \Omega). $$

We can see in [12], an explicit value of $ C\left (\dfrac{a}{b}\right ) =\sqrt {\dfrac{a}{b}} $. In his proof, Shafrir has used the Stokes formula and an isoperimetric inequality, see [6]. For Chen-Lin, they have used the blow-up analysis combined with some geometric type inequality for the integral curvature.

\bigskip

Now, if we suppose $ (V_i)_i $ uniformly Lipschitzian with $ A $ the
Lipschitz constant, then, $ C(a/b)=1 $ and $ c=c(a, b, A, K, \Omega)
$, see Brezis-Li-Shafrir [7]. This result was extended for
H\"olderian sequences $ (V_i)_i $ by Chen-Lin, see  [12]. Also, we
can see in [18], an extension of the Brezis-Li-Shafrir to compact
Riemann surface without boundary. We can see in [19] explicit form,
($ 8 \pi m, m\in {\mathbb N}^* $ exactly), for the numbers in front of
the Dirac masses, when the solutions blow-up. Here, the notion of isolated blow-up point is used.

\bigskip

In [8], Brezis and Merle proposed the following Problem:

\smallskip

{\bf Problem} {\it (Brezis-Merle [8])}.{\it If $ (u_i)_i $ and $ (V_i)_i $ are two sequences of functions relatively to the previous problem $ (P) $ with, 

$$ 0 \leq V_i \to V \,\, {\rm in } \,\,  C^0(\bar \Omega). $$

$$ \int_{\Omega} e^{u_i} dy  \leq C, $$

Is it possible to prove that:

$$ \sup_{\Omega} u_i \leq c = c(C, V, \Omega) \,\, ? $$}

Here, we assume more regularity on $ V_i $, we suppose that $ V_i \geq 0 $ is $ C^s $ ($ s$-holderian) $ 1/2 < s \leq 1)$ and when we have a boundary singularity. We give the answer where $ b C < 16 \pi $ for an equation with boundary singularity.

Our main results are:

\smallskip

\begin{Theorem} Assume $ \Omega=B_1(0) $, $ x_0 \in \partial \Omega $, $ \alpha  \in (0, 1/2) $, and,

$$ \int_{B_1(0)} |x-x_0|^{-2\alpha}V_i e^{u_i}  dy  \leq 16\pi-\epsilon, \,\, \epsilon>0,$$

$$ u_i(x_i)=\sup_{B_1(0)} u_i \to + \infty. $$

 There is a sequences $ (x_i^0)_i, (\delta_i^0) $, such that:

$$  (x_i^0)_i \equiv (x_i)_i, \,\, \delta_i^0=\delta_i=d(x_i, \partial B_1(0)) \to 0, $$

and,

$$ u_i(x_i)=\sup_{B_1(0)} u_i \to + \infty, $$

$$ u_i(x_i) + 2\log \delta_i- 2 \alpha \log d(x_i, x_0)  \to + \infty, $$

$$ \forall \,\, \epsilon >0, \,\,\, \limsup_{i \to + \infty} \int_{B(x_i, \delta_i\epsilon)} |x-x_0|^{-2\alpha}V_i e^{u_i}  dy \geq 4 \pi >0 .$$

 If we assume:

$$ V_i \to V\,\, {\rm in } \,\, C^0(\bar B_1(0)), $$

then,

$$ \forall \,\, \epsilon >0, \,\, \sup_{B_1(0)- B(x_i, \delta_i \epsilon)} u_i \leq C_{\epsilon} $$

$$  \forall \,\, \epsilon >0, \,\,\, \limsup_{i \to + \infty} \int_{B(x_i, \delta_i\epsilon)} |x-x_0|^{-2\alpha}V_i e^{u_i}  dy = 8 \pi . $$

And, thus, we have the following convergence in the sense of distributions:

$$ \int_{B_1(0)}  |x-x_0|^{-2\alpha}V_i e^{u_i}  dy  \to \int_{B_1(0)}  |x-x_0|^{-2\alpha}V e^{u}  dy + 8\pi \delta_{x_0}. $$

\end{Theorem}

\begin{Theorem}Assume that, $ V_i  $ is uniformly $ s- $holderian with $ 1/2 < s \leq 1$, $ x_0 \in \partial \Omega $, $ \alpha  \in (0, 1/2) $, and,

$$ \int_{B_1(0)} |x-x_0|^{-2\alpha}V_i e^{u_i}  dy  \leq 16\pi-\epsilon, \,\, \epsilon>0,$$

then we have:

$$ \sup_{\Omega} u_i \leq c=c(b, C, A, s,\alpha, \epsilon, x_0,  \Omega). $$

where $ A $ is the h\"olderian constant of $ V_i $.

\end{Theorem}

\section{Proofs of the results} 

\underbar {Proofs of the theorems:}

\bigskip

Without loss of generality, we can assume that $ \Omega= B_1(0) $ the unit ball centered on the origin.

\bigskip

Here, $ G $ is the Green function of the Laplacian with Dirichlet condition on $ B_1(0) $. We have (in complex notation):

$$ G(x, y)=\dfrac{1}{2 \pi} \log \dfrac{|1-\bar xy|}{|x-y|}, $$ 

Since $ u_i \in W^{1,1}_0(\Omega) $ and $ \alpha \in (0, 1/2) $, we have by the Brezis-Merle result and the elliptic estimates, (see [1]):

$$ u_i \in C^2(\Omega) \cap W^{2, p}(\Omega) \cap C^{1,\epsilon}(\bar \Omega) $$

for all $ 2 < p <+\infty $.

Set,

$$ v_i(x) =\int_{B_1(0)} G(x,y) V_i(y) |x-x_0|^{-2\alpha}e^{u_i(y)} dy . $$

We decompose $ v_i $ in two terms (Newtionian potential):

$$ v_i^1(x) =\int_{B_1(0)} -\dfrac{1}{2 \pi} \log |x-y| V_i(y) |x-x_0|^{-2\alpha}e^{u_i(y)} dy, $$

and,

$$ v_i^2(x) =\int_{B_1(0)} \dfrac{1}{2 \pi} \log |1-\bar x y| V_i(y) |x-x_0|^{-2\alpha}e^{u_i(y)} dy, $$

According to the proof in the book of Gilbarg-Trudinger see [15], $ v_i^1 $, $ v_i^2 $ and thus $ v_i $ are $ C^1(\bar \Omega) $. Indeed, we use the same proof as in [15] (Chapter 4, Newtonian potential), we have for the approximate function $ \partial v_{i, \epsilon} $ terms of type $ O(\epsilon^{1-2\alpha}\log \epsilon) + O(\epsilon^{1-2\alpha}) $. Since $\alpha < 1/2$, this term tends to $ 0 $.

We use this fact and the maximum principle to have $ v_i=u_i $.

Also, we can use integration by part (the Green representation formula, see its proof in the first chapter of [15]) to have in $ \Omega $ (and not $ \bar \Omega  $):

$$ u_i(x)=-\int_{B_1(0)} G(x,y)\Delta u_i(y) dy=\int_{B_1(0)} G(x,y) V_i(y) |x-x_0|^{-2\alpha}e^{u_i(y)} dy. $$
  
We write,

$$ u_i(x_i)=\int_{\Omega} G(x_i,y) |x-x_0|^{-2\alpha}V_i(y) e^{u_i(y)} dx=\int_{\Omega-B(x_i,\delta_i/2)} G(x_i,y) |x-x_0|^{-2\alpha}V_i e^{u_i(y)} dy + $$
 $$ + \int_{B(x_i,\delta_i/2)} G(x_i,y) |x-x_0|^{-2\alpha}V_ie^{u_i(y)} dy  $$ 
 
According to the maximum principle,  the harmonic function $ G(x_i, .) $ on $ \Omega-B(x_i,\delta_i/2) $ take its maximum on the boundary of $ B(x_i,\delta_i/2) $, we can compute this maximum:

$$ G(x_i, y_i)=\dfrac{1}{2 \pi} \log \dfrac{|1-\bar x_iy_i|}{|x_i-y_i|} =\dfrac{1}{2\pi} \log \dfrac{|1-\bar x_i(x_i+\delta_i\theta_i)|}{|\delta_i/2| } =\dfrac{1}{2\pi} \log (2|(1+|x_i|)+\theta_i|) < + \infty $$

with $ |\theta_i|= 1/2 $.

\bigskip

Thus,

$$ u_i(x_i) \leq C+  \int_{B(x_i,\delta_i/2)} G(x_i,y) |x-x_0|^{-2\alpha}V_ie^{u_i(y)} dy  \leq C + e^{u_i(x_i)-2 \alpha \log d(x_i, x_0)}  \int_{B(x_i,\delta_i/2)} G(x_i,y) dy $$

Now, we compute $ \int_{B(x_i,\delta_i/2)} G(x_i,y) dy  $

we set in polar coordinates,

$$ y=x_i+\delta_i t \theta  $$

we find:

$$ \int_{B(x_i,\delta_i/2)} G(x_i,y) dy =\int_{B(x_i,\delta_i/2)} \dfrac{1}{2 \pi} \log \dfrac{|1-\bar x_iy|}{|x_i-y|} =\dfrac{1}{2 \pi}\int_0^{2\pi}\int_0^{1/2} \delta_i^2 \log \dfrac{|1-\bar x_i(x_i+\delta_i t \theta)|}{\delta_i t} t dt d \theta = $$

$$ = \dfrac{1}{2 \pi}\int_0^{2\pi}\int_0^{1/2} \delta_i^2 (\log (|1+|x_i|+ t \theta|)-\log t) t dt d \theta \leq C \delta_i^2.$$

Thus,

$$ u_i(x_i) \leq C+ C \delta_i^2 e^{u_i(x_i)-2 \alpha \log d(x_i, x_0)} , $$

which we can write, because $ u_i(x_i) \to +\infty $,

$$ u_i(x_i) \leq C' \delta_i^2 e^{u_i(x_i)-2 \alpha \log d(x_i, x_0)}, $$

We can conclude that:

$$ u_i(x_i)+2 \log \delta_i-2 \alpha \log d(x_i, x_0) \to + \infty. $$

Since in $ B(x_i, \delta_i \epsilon) $, $ d(x, x_0) $ is equivalent to $ d(x_i,x_0) $ we can consider the following functions:

$$ v_i(y)=u_i(x_i+\delta_i y)+2\log \delta_i-2 \alpha \log d(x_i, x_0), \quad y \in B(0,1/2) $$

The function satisfies all conditions of the Brezis-Merle hypothesis, we can conclude that, on each compact set:

$$ v_i \to -\infty $$

we can assume, without loss of generality that for $ 1/2 > \epsilon >0 $, we have:

$$ v_i \to -\infty, \quad y \in  B(0, 2\epsilon)-B(0,\epsilon), $$

\begin{Lemma} 

For all $ 1/4> \epsilon >0 $, we have:

$$ \sup_{B(x_i, (3/2)\delta_i\epsilon)-B(x_i, \delta_i \epsilon)} u_i \leq C_{\epsilon}. $$

\end{Lemma}

\underbar{Proof of the lemma}

\bigskip

Let $ t'_i $ and $ t_i $ the points of $ B(x_i, 2\delta_i\epsilon)-B(x_i, (1/2)\delta_i \epsilon) $ and $ B(x_i, (3/2)\delta_i\epsilon)-B(x_i, \delta_i \epsilon) $ respectively where $ u_i $ takes its maximum.

According to the Brezis-Merle work, we have:

$$ u_i(t'_i)+2\log \delta_i-2 \alpha \log d(x_i, x_0) \to - \infty $$

We write,

$$ u_i(t_i)=\int_{\Omega} G(t_i,y) |x-x_0|^{-2\alpha}V_i(y) e^{u_i(y)} dx=\int_{\Omega-B(x_i,2\delta_i\epsilon)} G(t_i,y) |x-x_0|^{-2\alpha}V_i e^{u_i(y)} dy + $$

$$+ \int_{B(x_i,2\delta_i\epsilon)-B(x_i, (1/2)\delta_i\epsilon)} G(t_i,y) |x-x_0|^{-2\alpha}V_ie^{u_i(y)} dy+ $$

$$ +  \int_{B(x_i, (1/2)\delta_i\epsilon)} G(t_i,y) |x-x_0|^{-2\alpha}V_ie^{u_i(y)} dy $$ 

But, in the first and the third integrale, the point $ t_i $ is far from the singularity $ x_i $ and we know that the Green function is bounded. For the second integrale, after a change of variable, we can see that this integale is bounded by (we take the supremum in the annulus and use Brezis-Merle theorem)

$$ \delta_i^2 e^{u_i(t'_i)-2 \alpha \log d(x_i, x_0)} \times  I_j $$

where $ I_j $ is a Jensen integrale (of the form $ \int_0^{1}  \int_0^{2\pi}  (\log (|1+|x_i|+ t \theta)-\log |\theta_i-t\theta|) t dt d \theta $ which is bounded ). 

we conclude the lemma.

From the lemma, we see that far from the singularity the sequence is bounded, thus if we take the supremum on the set $ B_1(0)-B(x_i, \delta_i \epsilon) $ we can see that this supremum  is bounded and thus the sequence of functions is uniformly bounded or tends to infinity and we use the same arguments as for $ x_i $  to conclude that around this point and far from the singularity, the seqence is bounded.

The process will be finished , because, according to Brezis-Merle estimate, around each supremum constructed and tending to infinity, we have:

$$ \forall \,\, \epsilon >0, \,\,\, \limsup_{i \to + \infty} \int_{B(x_i, \delta_i\epsilon)} |x-x_0|^{-2\alpha}V_i e^{u_i}  dy \geq 4 \pi >0 .$$

Finaly, with this construction, we have a finite number of "exterior "blow-up points and outside the singularities the sequence is bounded uniformly, for example, in the case of one "exterior" blow-up point, we have:

$$ u_i(x_i) \to + \infty  $$

$$ \forall \,\, \epsilon >0, \,\, \sup_{B_1(0)-B(x_i, \delta_i \epsilon)} u_i \leq C_{\epsilon} $$

$$ \forall \,\, \epsilon >0, \,\,\, \limsup_{i \to + \infty} \int_{B(x_i, \delta_i\epsilon)} |x-x_0|^{-2\alpha}V_i e^{u_i}  dy \geq 4 \pi >0 .$$

$$ x_i \to x_0 \in \partial B_1(0). $$

{\bf Remark}: For the general case, the process of quantization can be extended to more than one blow-up points.

We have the following lemma:

\begin{Lemma} Each $ \delta_i^k $ is of order $ d(x_i^k, \partial B_1(0)) $. Namely: there is a positive constant $ C >0 $ such that for $ \epsilon >0 $ small enough:

$$ \delta_i^k \leq d(x_i^k, \partial B_1(0)) \leq (2+\dfrac{C}{\epsilon})\delta_i^k. $$

\end{Lemma}

\underbar{Proof of the lemma}

Now, if we suppose that there is another "exterior" blow-up $ (t_i)_i $, we have, because $ (u_i)_i $ is uniformly bounded in a neighborhood of $ \partial B(x_i, \delta_i\epsilon) $, we have :

$$ d(t_i,\partial B(x_i, \delta_i\epsilon)) \geq \delta_i\epsilon $$

If we set,

$$ \delta'_i= d(t_i, \partial (B_1(0)-B(x_i, \delta_i\epsilon)))= \inf \{d(t_i,\partial B(x_i, \delta_i\epsilon)), d(t_i,\partial (B_1(0))) \} $$

then, $ \delta'_i  $ is of order $ d(t_i, \partial B_1(0)) $. To see this, we write:

$$ d(t_i,\partial B_1(0)) \leq d(t_i,\partial B(x_i, \delta_i\epsilon))+d( \partial B(x_i, \delta_i\epsilon), x_i) + d(x_i, \partial B_1(0)), $$

Thus,

$$ \dfrac{d(t_i,\partial B_1(0))}{d(t_i,\partial B(x_i, \delta_i\epsilon))} \leq 2+\dfrac{1}{\epsilon}, $$

Thus,

$$   \delta'_i \leq d(t_i,\partial B_1(0)) \leq \delta'_i (2+\dfrac{1}{\epsilon}). $$

\bigskip

Now, the general case follow by induction. We use the same argument for three, four,..., $ n $ blow-up points. 

We have, by induction and, here we use the fact that $ u_i $ is uniformly bounded outside a small ball centered at $ x_i^j, j=0,\ldots,k-1 $:

$$ \delta_i^j  \leq d(x_i^j, \partial B_1(0))\leq C_1 \delta_i^j, \,\, j=0,\ldots, k-1,  $$.

$$ d(x_i^k, \partial B(x_i^j, \delta_i^j \epsilon/2))\geq \epsilon \delta_i^j, \epsilon >0,\,\, j=0,\ldots, k-1,  $$.

and let's consider $ x_i^k $ such that:

$$ u_i(x_i^k)=\sup_{B_1(0)- \cup_{j=0}^{k-1}  B(x_i^j, \delta_i^j \epsilon)} u_i \to + \infty, $$

take,

$$ \delta_i^k= \inf \{d(x_i^k, \partial B_1(0)), d(x_i^k, \partial (B_1(0)-\cup_{j=0}^{k-1}  B(x_i^j, \delta_i^j \epsilon/2)) \}, $$

if, we have,

$$ \delta_i^k= d(x_i^k, \partial B(x_i^j, \delta_i^j \epsilon/2)), \,\, j\in\{0, \ldots,  k-1 \}.  $$ 

Then,

$$ \delta_i^k \leq d(x_i^k, \partial B_1(0)) \leq $$

$$ \leq d(x_i^k, \partial B(x_i^j, \delta_i^j \epsilon/2))+ d(\partial B(x_i^j, \delta_i^j \epsilon/2), x_i^j)+ d(x_i^j, \partial B_1(0))  $$

$$ \leq (2+\dfrac{C_1}{\epsilon})\delta_i^k. $$

To apply lemma 2.1 for $ m $ blow-up points, we use an induction:

\bigskip

We do directly the same approch for $ t_i $ as $ x_i $ by using directly the Green function of the unit ball.

\bigskip

If we look to the blow-up points, we can see, with this work that, after finite steps, the sequence will be bounded outside a finite number of balls , because of Brezis-Merle estimate:

$$  \forall \,\, \epsilon >0, \,\,\, \limsup_{i \to + \infty} \int_{B(x_i^k, \delta_i^k\epsilon)} |x-x_0|^{-2\alpha}V_i e^{u_i}  dy \geq 4 \pi >0 .$$

Here, we can take the functions:

$$ u_i^k(y)=u_i(x_i^k+\delta_i^k y)+2\log \delta_i^k-2\alpha \log d(x_i,x_0). $$

Indeed, by corollary 4 of the paper of Brezis-Merle, if we have:

$$ \limsup_{i \to + \infty} \int_{B(x_i^k, \delta_i^k\epsilon)} |x-x_0|^{-2\alpha}V_i e^{u_i}  dy \leq 4 \pi-\epsilon_0 <4\pi, $$

then, $ (u_i^k)^+ $ would be bounded and this contradict the fact that $ u_i^k(0)\to +\infty $.

\smallskip

Finaly, we can say that, there is a finite number of sequences $ (x_i^k)_i, (\delta_i^k) , 0 \leq k \leq m $, such that:

$$  (x_i^0)_i \equiv (x_i)_i, \,\, \delta_i^0=\delta_i=d(x_i, \partial B_1(0)), $$

$$ (x_i^1)_i \equiv (t_i)_i, \,\, \delta_i^1=\delta'_i=d(t_i, \partial (B_1(0)-B(x_i, \delta_i \epsilon)), $$

and each $\delta_i^k $ is of order $ d(x_i^k, \partial B_1(0)) $.

and,

$$ u_i(x_i^k)=\sup_{B_1(0)- \cup_{j=0}^{k-1}  B(x_i^j, \delta_i^j \epsilon)} u_i \to + \infty, $$

$$ u_i(x_i^k) + 2\log \delta_i^k-2 \alpha \log d(x_i^k, x_0) \to + \infty, $$

$$ \forall \,\, \epsilon >0, \,\, \sup_{B_1(0)- \cup_{j=0}^{m}  B(x_i^j, \delta_i^j \epsilon)} u_i \leq C_{\epsilon} $$

$$ \forall \,\, \epsilon >0, \,\,\, \limsup_{i \to + \infty} \int_{B(x_i^k, \delta_i^k\epsilon)} |x-x_0|^{-2\alpha}V_i e^{u_i}  dy \geq 4 \pi >0 .$$

\underbar {The work of YY.Li-I.Shafrir}

\bigskip

Since in $ B(x_i, \delta_i \epsilon) $, $ d(x, x_0) $ is equivalent to $ d(x_i,x_0) $ we can consider the following functions:

 $$ v_i(y)=u_i(x_i+\delta_i y)+2\log \delta_i-2 \alpha \log d(x_i, x_0). $$

With the previous method, we have a finite number of "exterior" blow-up points (perhaps the same) and the sequences tend to the boundary. With the aid of proposition 1 of the paper of Li-Shafrir, we see that around each exterior blow-up, we have a finite number of "interior" blow-ups. Around, each exterior blow-up, we have after rescaling with $ \delta_i^k $, the same situation as around a fixed ball with positive radius.  If we assume:

$$ V_i \to V\,\, {\rm in } \,\, C^0(\bar B_1(0)), $$

then,

$$  \forall \,\, \epsilon >0, \,\,\, \limsup_{i \to + \infty} \int_{B(x_i^k, \delta_i^k\epsilon)}|x-x_0|^{-2\alpha} V_i e^{u_i}  dy = 8 \pi m_k, \,\, m_k \in  {\mathbb N}^*.$$

And, thus, we have the following convergence in the sense of distributions:

$$ \int_{B_1(0)} |x-x_0|^{-2\alpha}V_i e^{u_i}  dy  \to \int_{B_1(0)} |x-x_0|^{-2\alpha}V e^{u}  dy +\sum_{k=0}^m 8 \pi m_k' \delta_{x_0^k} ,\, m_k' \in  {\mathbb N}^*, \,\, x_0^k \in \partial B_1(0). $$

\underbar {Consequence: using a Pohozaev-type identity, proof of theorem 2}

\bigskip

By a conformal transformation, we can assume that our domain $ \Omega = B^+ $ is a half ball centered at the origin, $ B^+=\{x, |x|\leq 1,  x_1 \geq 0 \} $, and, $ x_0=0 $. In this case the normal at the boundary is $ \nu=(-1,0) $ and $ u_i(0, x_2) \equiv 0 $. Also, we set $ x_i $ the blow-up point and $ x_i^2=(0,x_i^2) $ and $ x_i^1=(x_i^1,0) $ respectevely the second and the first part of $ x_i $. Let $ \partial B^+ $ the part of the boundary for which $ u_i $ and its derivatives are uniformly bounded and thus converge to the corresponding function.

\bigskip

\underbar {The case of one blow-up point:}

\bigskip

\begin{Theorem} If  $ V_i $ is s-Holderian with $ 1/2 < s  \leq 1$ and,

$$ \int_{\Omega} |x|^{-2\alpha}V_i e^{u_i} dy \leq 16\pi-\epsilon , \,\, \epsilon >0, $$

we have :
 
 $$ 2(1-2\alpha)V_i(x_i)\int_{B(x_i, \delta_i \epsilon)} |x|^{-2\alpha}e^{u_i} dy=o(1), $$
 
 which means that there is no blow-up points.

\end{Theorem}

\bigskip

\underbar{Proof of the theorem}

\bigskip

In order to use the Pohozaev identity we need to have a good function $ u_i $, since $ \alpha \in (0,1/2) $, we have a function $ u_i $ such that:

$$ u_i \in C^2(\Omega) \cap W^{2, p}(\Omega) \cap C^1(\bar \Omega) $$

Thus,

$$ \partial_j u_i, \partial_k u_i \in W^{1, p}(\Omega) \cap C^0(\bar \Omega). $$

Thus, we can use integration by parts, in fact we have for the following product (here $ "." $ is the usual product of function):

$$ \partial_j u_i . \partial_k u_i \in W^{1, p}(\Omega) \cap C^0(\bar \Omega). $$

The Pohozaev identity gives us the following formula:

$$ \int_{\Omega} <(x-x_2^i) |\nabla u_i>(-\Delta u_i) dy =\int_{\Omega} <(x-x_2^i) |\nabla u_i> |x|^{-2\alpha}V_i e^{u_i} dy = A_i $$

$$ A_i= \int_{\partial B^+} <(x-x_2^i) |\nabla u_i> <\nu |\nabla u_i> d\sigma + \int_{\partial B^+} <(x-x_2^i) |\nu>|\nabla u_i|^2d\sigma $$

We can write it as:

$$ \int_{\Omega} <(x-x_2^i) |\nabla u_i> (V_i-V_i(x_i))|x|^{-2\alpha} e^{u_i} dy =A_i+  V_i(x_i) \int_{\Omega} <(x-x_2^i) |\nabla u_i>|x|^{-2\alpha} e^{u_i} dy = $$

$$ = A_i+ V_i(x_i)\int_{\Omega} <(x-x_2^i)|x|^{-2\alpha} |\nabla (e^{u_i})> dy  $$

And, if we integrate by part the second term, we have (because $ x_1=0 $ on the boundary and $ \nu_2=0 $):

$$ \int_{\Omega} <(x-x_2^i) |\nabla u_i> (V_i-V_i(x_i))|x|^{-2\alpha} e^{u_i} dy =-2(1-\alpha)V_i(x_i)\int_{\Omega} |x|^{-2\alpha}e^{u_i} dy + $$

$$ +2\alpha V_i(x_i)\int_{B(x_i, \delta_i \epsilon)} x_2 x_2^i |x|^{-2\alpha-2}e^{u_i} dy+2\alpha V_i(x_i)\int_{\Omega - B(x_i, \delta_i \epsilon)} x_2 x_2^i |x|^{-2\alpha-2}e^{u_i} dy + B_i $$

where $ B_i $ is,

$$ B_i=V_i(x_i)\int_{\partial B^+}  <(x-x_2^i) |\nu>|x|^{-2\alpha} e^{u_i} dy $$

applying the same procedure to $ u $, we can write:

$$ \int_{\Omega} <(x-x_2^i) |\nabla u> (V-V(0))|x|^{-2\alpha} e^{u} dy =-2(1-\alpha)V(0)\int_{\Omega} |x|^{-2\alpha}e^{u} dy + $$

$$ +2\alpha V(0)\int_{B(x_i, \delta_i \epsilon)} x_2 x_2^i |x|^{-2\alpha-2}e^{u} dy+2\alpha V(0)\int_{\Omega - B(x_i, \delta_i \epsilon)} x_2 x_2^i |x|^{-2\alpha-2}e^{u} dy + B, $$

with,

$$ B=V(0)\int_{\partial B^+}  <(x-x_2^i) |\nu> |x|^{-2\alpha}e^{u} dy $$

we use the fact that, $ u_i $ is bounded outside $ B(x_i, \delta_i \epsilon) $ and the convergence of $ u_i $ to $ u $ on compact set of $ \bar \Omega-\{0\} $, and the fact that $ \alpha \in (0,1/2) $, to write the following:

$$ 2(1-2\alpha)V_i(x_i)\int_{B(x_i, \delta_i \epsilon)} |x|^{-2\alpha}e^{u_i} dy+o(1)= $$
$$ = \int_{\Omega} <(x-x_2^i) |\nabla u_i> (V_i-V_i(x_i)) |x|^{-2\alpha}e^{u_i} dy-\int_{\Omega} <(x-x_2^i) |\nabla u> (V-V(0)) |x|^{-2\alpha}e^{u} dy+ $$

$$ + (A_i-A)+ (B_i-B), $$

where $ A $ and $ B $ are,

$$ A= \int_{\partial B^+} <(x-x_2^i) |\nabla u> <\nu |\nabla u> d\sigma + \int_{\partial B^+} <(x-x_2^i) |\nu>|\nabla u|^2d\sigma $$

$$ B=V(0)\int_{\partial B^+}  <(x-x_2^i) |\nu> |x|^{-2\alpha}e^{u} dy $$

and, because of the uniform convergence of $ u_i $ and its derivatives on $ \partial B^+ $, we have:

$$ A_i-A = o(1) \,\,\, {\rm and}\,\,\, B_i-B = o(1) $$

which we can write as:

$$  2(1-2\alpha)V_i(x_i)\int_{B(x_i, \delta_i \epsilon)} |x|^{-2\alpha}e^{u_i} dy+o(1)= $$

$$ = \int_{\Omega} <(x-x_2^i) |\nabla (u_i-u)> (V_i-V_i(x_i)) |x|^{-2\alpha}e^{u_i} dy +  $$

$$ +  \int_{\Omega} <(x-x_2^i) |\nabla u> (V_i-V_i(x_i)) |x|^{-2\alpha}(e^{u_i}-e^u) dy + $$

$$  + \int_{\Omega} <(x-x_2^i) |\nabla u> (V_i-V_i(x_i)- (V-V(0)))|x|^{-2\alpha} e^{u} dy + o(1) $$ 

We can write the second term as:

 $$ \int_{\Omega} <(x-x_2^i) |\nabla u> (V_i-V_i(x_i)) |x|^{-2\alpha}(e^{u_i}-e^u) dy = \int_{\Omega-B(0,\epsilon) } <(x-x_2^i) |\nabla u> (V_i-V_i(x_i)) (e^{u_i}-e^u)|x|^{-2\alpha} dy + $$
 
 $$ + \int_{B(0,\epsilon)} <(x-x_2^i) |\nabla u> (V_i-V_i(x_i)) (e^{u_i}-e^u)|x|^{-2\alpha} dy = o(1), $$

because of the uniform convergence of $ u_i $ to $ u $ outside a region which contain the blow-up and the uniform convergence of $ V_i $. For the third integral we have the same result:

$$ \int_{\Omega} <(x-x_2^i) |\nabla u> (V_i-V_i(x_i)- (V-V(0)))|x|^{-2\alpha} e^{u} dy =o(1), $$

because of the uniform convergence of $ V_i  $ to $ V $.

Now, we look to the first integral:

$$ \int_{\Omega} <(x-x_2^i) |\nabla (u_i-u)> (V_i-V_i(x_i)) |x|^{-2\alpha}e^{u_i} dy,  $$

we can write it as:

$$ \int_{\Omega} <(x-x_2^i) |\nabla (u_i-u)> (V_i-V_i(x_i))|x|^{-2\alpha} e^{u_i} dy =\int_{\Omega} <(x-x_i) |\nabla (u_i-u)> (V_i-V_i(x_i)) |x|^{-2\alpha}e^{u_i} dy + $$

$$ + \int_{\Omega} <x_1^i |\nabla (u_i-u)> (V_i-V_i(x_i))|x|^{-2\alpha} e^{u_i} dy, $$

Thus, we have proved by using the Pohozaev identity  the following equality:

$$ \int_{\Omega} <(x-x_i) |\nabla (u_i-u)> (V_i-V_i(x_i)) |x|^{-2\alpha}e^{u_i} dy + $$

$$ + \int_{\Omega} <x_1^i |\nabla (u_i-u)> (V_i-V_i(x_i)) |x|^{-2\alpha}e^{u_i} dy = $$

$$ = 2(1-2\alpha)V_i(x_i)\int_{B(x_i, \delta_i \epsilon)} |x|^{-2\alpha}e^{u_i} dy+o(1)$$

We can see, because of the uniform boundedness of $ u_i $ outside $ B(x_i, \delta_i \epsilon) $ and the fact that :

$$ ||\nabla (u_i-u) ||_1 = o(1), $$

it is sufficient to look to the integral on $ B(x_i, \delta_i \epsilon) $.

\bigskip

Assume that we are in the case of one blow-up, it must be $ (x_i) $ and isolated, we can write the following inequality as a consequence of YY.Li-I.Shafrir result:

$$ u_i(x)+ 2\log |x-x_i|-2\alpha \log d(x,0) \leq C,  $$

We use this fact and the fact that $ V_i $ is s-holderian to have that, on $ B(x_i, \delta_i \epsilon) $,

$$ |(x-x_i)(V_i-V_i(x_i)) |x|^{-2 \alpha}e^{u_i}| \leq \dfrac{C}{ |x-x_i|^{1-s}} \in L^{(2-\epsilon')/(1-s)}, \,\, \forall \,\, \epsilon' >0, $$

and, we use the fact that:

$$ ||\nabla (u_i-u) ||_q = o(1), \,\, \forall \,\, 1 \leq q < 2$$

to conclude by the Holder inequality that:

$$ \int_{B(x_i, \delta_i \epsilon)} <(x-x_2^i) |\nabla (u_i-u)> (V_i-V_i(x_i)) |x|^{-2 \alpha}e^{u_i} dy = o(1), $$

For the other integral, namely:

$$  \int_{B(x_i, \delta_i \epsilon)} <x_1^i |\nabla (u_i-u)> (V_i-V_i(x_i)) |x|^{-2 \alpha}e^{u_i} dy , $$

We use the fact that, because our domain is a half ball, and the $ \sup+\inf  $ inequality to have:

$$ x_1^i= \delta_i, $$

$$ u_i(x)+4 \log \delta_i - 4\alpha \log d(x,0) \leq C $$

and,

$$ |x|^{-s \alpha} e^{ (s/2)u_i(x)} \leq |x-x_i|^{-s},  $$

$$ |V_i-V_i(x_i)|\leq |x-x_i|^s, $$

Finaly, we have:

$$  | \int_{B(x_i, \delta_i \epsilon)} <x_1^i |\nabla (u_i-u)> (V_i-V_i(x_i)) |x|^{-2 \alpha}e^{u_i} dy |\leq $$

$$ \leq C \int_{B(x_i, \delta_i \epsilon)}  |\nabla (u_i-u)|(|x|^{-2 \alpha}e^{u_i})^{(3/4-s/2)}, $$

But in the second member, for $ 1/2 < s \leq 1 $, we have $ q_s=1/(3/4-s/2) >2 $ and thus $ q_s' <2 $ and,

$$ (|x|^{-2\alpha }e^{u_i})^{3/4-s/2} \in L^{q_s} $$

$$ |x|^{-2\alpha (3/4-s/2)}e^{((3/4)-(s/2))u_i} \in L^{q_s} $$

$$  ||\nabla (u_i-u) ||_{q_s'}  = o(1), \,\, \forall \,\, 1 \leq q_s' < 2, $$

one conclude that:

$$  \int_{B(x_i, \delta_i \epsilon)} <x_1^i |\nabla (u_i-u)> (V_i-V_i(x_i)) |x|^{-2\alpha}e^{u_i} dy =o(1) $$

 Finaly, with this method, we conclude that, in the case of one blow-up point and $ V_i $ is s-Holderian with $ 1/2 < s  \leq 1$ :

 $$  2(1-2\alpha)V_i(x_i)\int_{B(x_i, \delta_i \epsilon)} |x|^{-2\alpha}e^{u_i} dy=o(1)$$
 
 which means that there is no blow-up, which is a contradiction.
 
\bigskip 
 
 Finaly, for one blow-up point and $ V_i  $ is  is s-Holderian with $ 1/2 < s  \leq 1$, the sequence $ (u_i) $ is uniformly bounded on $ \Omega $.

\end{document}